\documentclass[12pt]{article}
\usepackage{latexsym,amssymb,amsmath,amsfonts,amsthm,mathrsfs,cite}

\newtheorem{thm}{Theorem}[section]
\newtheorem{prop}[thm]{Proposition}
\newtheorem{lem}[thm]{Lemma}
\newtheorem{cor}[thm]{Corollary}
\def\pf{\noindent{\it Proof.} }
\newcommand{\E}{\mathrm{E}}
\newcommand{\V}{\mathrm{V}}
\def\qed{\nopagebreak\hfill{\rule{4pt}{7pt}}\medbreak}

\makeatletter \@addtoreset{equation}{section} \makeatother

\begin{document}
\begin{center}
{\Large\bf
The Limiting Distribution of the Number of

Block Pairs in Type $B$ Set Partitions}
\end{center}

\begin{center}
David G. L. Wang

Beijing International Center for Mathematical Research\\
Peking University, Beijing 100871, P. R. China\\
wgl@math.pku.edu.cn
\end{center}

\begin{abstract}
It is a classical result of Harper that the limiting distribution of the
number of blocks in partitions of the set $\{1, 2,\ldots, n\}$ is normal.
In this paper, using the saddle point method
we prove the normality of the limiting distribution
of the number of block pairs in set partitions of type $B_n$.
Moreover, we obtain that the limiting
distribution of the number of block pairs
in $B_n$-partitions without zero-block is also normal.
\end{abstract}

\section{Introduction}

This paper is concerned with the limiting distribution
of the number of block pairs of type $B_n$ set partitions.
For ordinary set partitions, Harper~\cite{Har67} has established
the normality of the limiting distribution
of the number
of blocks in partitions of the set $\{1,2,\ldots,n\}$.
For the asymptotic behavior concerning with
the ordinary set partitions,
see~\cite{DP83,GS92,GS94,KOPRSSW99,OR85-1,Pittel97}.
For the study on the limiting distribution of other
combinatorial objects,
see Flajolet and Sedgewick's book~\cite{FS09}
for instance.

The lattice of ordinary set partitions can be regarded as the
intersection lattice for the hyperplane
arrangement corresponding
to the root system of type $A$,
see Bj\"orner and Brenti~\cite{BB05}
or Humphreys~\cite{Hum92}.
From this point of view,
type $B$ set partitions are a
generalization of ordinary partitions, see
Reiner~\cite{Rei97}.
To be more precise,
ordinary set partitions
encode the intersections of hyperplanes in the hyperplane
arrangement for the type $A$ root system,
 while the intersections
of subsets of hyperplanes from the type~$B$ hyperplane
arrangement can be encoded by type~$B$ set partitions,
see Bj\"{o}rner and Wachs~\cite{BW04}.
A type~$B_n$ set partition is a partition~$\pi$ of the set
\[
\{1, 2, \ldots, n, -1, -2, \ldots, -n\}
\]
such that for any block~$B$ of~$\pi$, $-B$~is also a
block of~$\pi$, and there is at most one block,
called zero-block, satisfying $B=-B$.
We call $(B,-B)$ a block pair of~$\pi$
if~$B$ is not a zero-block.

Let $M_{n,\,k}$ be the number of $B_n$-partitions
with $k$ block pairs.
It is easy to deduce the recurrence relation
\begin{equation}\label{RecM}
M_{n,\,k}=M_{n-1,\,k-1}+(2k+1)M_{n-1,\,k}.
\end{equation}
The main result of this paper is to derive the limiting distribution of
the number of block pairs in $B_n$-partitions based on the above recurrence formula.
Let~$\xi_n$ be the random variable of the number of
block pairs in $B_n$-partitions.
We shall prove that the limiting distribution of~$\xi_n$
is normal by using the saddle point method,
which was introduced by Schr\"odinger~\cite{Sch44},
see also~\cite{Odl95,GK82,deBru58,Ben74}.

This paper is organized as follows.
In Section 2, we present some facts about the saddle point of the
generating function of the number of $B_n$-partitions.
Section~3 is devoted to
deduce the normality of the limiting distribution of $\xi_n$.
Using the same technique,
we obtain the normality of the limiting distribution of
the number of block-pairs
in $B_n$-partitions without zero-block.

\section{Preliminary lemmas}

Let $M_n$ be the number of $B_n$-partitions.
In this section,
we give some lemmas which will be used
to derive an approximate formula for $M_n$.

Let $N_{n,\,k}$ be the number of $B_n$-partitions with $k$ block pairs but with no zero-block.
Denote by $N_n$ the number of $B_n$-partitions without zero-block.
It is easy to see that
\[
N_{n,\,k}=2^{n-k}S(n,k).
\]
Since
\begin{equation}
\sum_{n}S(n,k){x^n\over n!}={1\over k!}\left(e^x-1\right)^k,
\end{equation}see Stanley~\cite[page 34]{Sta97},
we find that
\begin{equation}\label{EgfN}
F_N(z)=\sum_{n\ge0}N_n\,{z^n\over n!}
=\exp\Bigl(\frac{e^{2z}-1}{2}\Bigr).
\end{equation}
It is also easy to see that
\[
M_n=\sum_{k}{n\choose k}N_k.
\]
It follows from~\eqref{EgfN} that
\begin{equation}\label{EgfM}
F_M(z)=\sum_{n\ge0}M_n\,{z^n\over n!}
=\exp\Bigl(\frac{e^{2z}-1}{2}+z\Bigr).
\end{equation}

The saddle point of $F_M(z)$ is defined to be
the value $z$ that minimizes $z^{-n}F_M(z)$,
i.e., the unique positive solution $r_1$ of the equation
\begin{equation}\label{r1}
r_1\bigl(e^{2r_1}+1\bigr)=n.
\end{equation}
Similarly,
the saddle point of $F_N(z)$
is the unique positive solution $r_0$ of the equation
\begin{equation}\label{r0}
r_0e^{2r_0}=n.
\end{equation}
For convenience, we consider the equation
\begin{equation}\label{r}
r\left(e^{2r}+c\right)=n.
\end{equation}
It reduces to~\eqref{r1} when $c=1$, and to~\eqref{r0} when $c=0$.
It is easy to deduce the following approximation for the unique positive
solution $r$ of~\eqref{r}.

\begin{lem}\label{lem_r}
Let $c$ be a nonnegative integer.
Let $r>0$ be the unique positive solution of equation~\eqref{r}. Then we have
\begin{align*}
r&={\log{n}\over 2}\biggl(1+O\Bigl({\log\log{n}\over\log{n}}\Bigr)\biggr),\\[5pt]
e^{2r}&=
{2n\over\log{n}}\biggl(1+O\Bigl({\log\log{n}\over\log{n}}\Bigr)\biggr).
\end{align*}
\end{lem}

We will also need the following lemma.

\begin{lem}\label{lem_tech}
Let $h(x)$ be a continuous function defined on the closed
interval $[a,b]$. Suppose that $h''(x)$ exists
in the open interval $(a,b)$. Then for any $c\in(a,b)$,
there exists $s\in(a,b)$ such that
\begin{equation}\label{eq}
{h(a)\over(a-b)(a-c)}+{h(b)\over(b-a)(b-c)}
+{h(c)\over(c-a)(c-b)}={h''(s)\over 2}.
\end{equation}
\end{lem}

\pf
Let
\begin{align*}
f_1(x)&=(a-b)h(x)+(b-x)h(a)+(x-a)h(b),\\[5pt]
g_1(x)&=(a-b)(b-x)(x-a).
\end{align*}
Then the left hand side of~\eqref{eq}
becomes $f_1(c)/g_1(c)$. Note that $f_1(a)=g_1(a)=0$.
By Cauchy's mean value theorem, there exists $s_1\in(a,c)$ such that
\[
{f_1(c)\over g_1(c)}={f_1'(s_1)\over g_1'(s_1)}
={f_2(a)-f_2(b)\over g_2(a)-g_2(b)},
\]
where $f_2(x)=h'(s_1)x-h(x)$ and $g_2(x)=x^2-2s_1x$.
Again, by Cauchy's mean value theorem,
there exist $s_2\in(a,b)$ and $s\in(a,b)$ such that
\[
{f_1(c)\over g_1(c)}={f_2'(s_2)\over g_2'(s_2)}
={h'(s_1)-h'(s_2)\over 2s_1-2s_2}={h''(s)\over 2}.
\]
This completes the proof. \qed

\begin{lem}\label{lem1}
Let~$c$ be a nonnegative integer, and $f(x)=x\bigl(e^{2x}+c\bigr)$.
Suppose that $t_i$ ($i=0,1,2$) is the unique positive number such that
$f(t_i)=n+i$.
Then we have
\begin{align}
&t_1-t_0={1\over 2n}-{1\over 4nt_0}+O\Bigl({1\over n\log^2{n}}\Bigr),
\label{61}\\[5pt]
&t_2-t_1={1\over 2n}-{1\over 4nt_0}+O\Bigl({1\over n\log^2{n}}\Bigr),
\label{62}\\[5pt]
&2t_1-t_0-t_2={1\over n^2}+O\Bigl({1\over n^2\log{n}}\Bigr),\notag\\[5pt]
&{1\over t_0}+{1\over t_2}-{2\over t_1}=O\Bigl({1\over n^2\log^2{n}}\Bigr).\notag
\end{align}
\end{lem}

\pf We first consider~\eqref{61}.
By Cauchy's mean value theorem,
there exists $t$ such that
$t_0<t<t_1$ and
\[
f(t_1)-f(t_0)=(t_1-t_0)f'(t).
\]
Since $f(t_1)-f(t_0)=1$ and $f'(t)=(2t+1)e^{2t}+c$, we have
\begin{equation}\label{56}
\frac{1}{(2t_1+1)e^{2t_1}+c}\le t_1-t_0\le\frac{1}{(2t_0+1)e^{2t_0}+c}.
\end{equation}
It can be seen that both ${1 \over (2t_1+1)e^{2t_1}+c}$ and $\frac{1}{(2t_0+1)e^{2t_0}+c}$
have the same estimate
\[
{1\over 2n}-{1\over 4nt_0}+O\Bigl({1\over n\log^2{n}}\Bigr).
\]
It follows that $t_1-t_0$ also has the above estimate.
Similarly, one can prove~\eqref{62}.
The last two approximations are consequences of Lemma~\ref{lem_tech}.
\qed

\section{The limiting distribution}

Recall that $\xi_n$ is the random variable
of the number of block pairs in a $B_n$-partition.
Denote by $\E(\xi_n)$ the expectation of $\xi_n$,
and $\V(\xi_n)$ the variance of $\xi_n$.
Below is the main result of this paper.

\begin{thm}\label{thm_DistM}
The limiting distribution of the random variable $\xi_n$
is normal. In other words, the random variable
\[
{\xi_n-\E(\xi_n)\over \sqrt{\V(\xi_n)}}
\]
has an asymptotically standard normal distribution as~$n$~tends to
infinity.
\end{thm}

There are various sufficient conditions on a random variable
which ensures a normal limiting distribution,
see Sachkov~\cite{Sac97}.
Let $\eta_n$ be a random variable of certain statistic of
some combinatorial objects on a set $A_n$.
Let $a_n(k)$ be the number of elements of $A_n$ with the statistic equal to $k$.
Consider the polynomial
\[
P_n(x)=\sum_ka_n(k)x^k.
\]
The following criterion
was used by Harper~\cite{Har67},
see also Bender~\cite{Ben73}.

\begin{prop}\label{prop_1}
The limiting distribution
of $\eta_n$ is normal, if the $P_n(x)$ distinct real roots
and the variance of $a_n(k)$ tends to infinity as $n\to\infty$.
\end{prop}

We shall prove Theorem~\ref{thm_DistM} with the aid of Proposition~\ref{prop_1}.
Recall that $M_{n,\,k}$ is the number of $B_n$-partitions
with $k$ block pairs. Consider the polynomial
\begin{equation}\label{eq_Mnx}
M_n(x)=\sum_{k}M_{n,\,k}\,x^k.
\end{equation}
For example,
\begin{align*}
M_1(x)&=1+x,\\[5pt]
M_2(x)&=1+4x+x^2, \\[5pt]
M_3(x)&=1+13x+9x^2+x^3.
\end{align*}

\begin{thm}\label{thm_RootM} For any $n\geq 1$,
the polynomial $M_n(x)$ has $n$ distinct real roots.
\end{thm}

\pf The proof is similar to the
proof of Harper for ordinary partitions.
We prove it by induction on $n$.
It is clear that the theorem holds for $n=1, 2$.
We assume that it holds for all $n\le m-1$, where $m\geq 3$.
Let
\begin{equation}\label{G}
G_n(x)=\sqrt{x}\,e^{x\over 2}M_n(x).
\end{equation}
Differentiating $M_n(x)$ with respect to $x$ and using the recurrence~\eqref{RecM},
we obtain that
\begin{equation}\label{Mm}
M_n(x)=(1+x)M_{n-1}(x)+2xM_{n-1}'(x).
\end{equation}
Multiplying both sides of~\eqref{Mm} by $\sqrt{x}\,e^{x\over 2}$ yields
\begin{equation}\label{Gm}
G_n(x)=2xG_{n-1}'(x).
\end{equation}
By the induction hypothesis, we may assume that
$M_{m-1}(x)$ has roots $x_1,x_2,\ldots,x_{m-1}$
where $x_1<x_2<\cdots<x_{m-1}<0$.
Observe that
\[
\lim\limits_{x\to-\infty}G_n(x)=0.
\]
From~\eqref{G} it can be seen that $G_{m-1}(x)$ has $m+1$ roots
\[
 -\infty, \ x_1,\ x_2,\ \ldots,\ x_{m-1},\ 0.
\]
By Rolle's theorem, in each of the $m$ open intervals
\[
(-\infty,x_1),\quad (x_1,x_2),\quad \ldots,\quad (x_{m-1},0),
\]
there exists a point $y$ such that $G_{m-1}'(y)=0$.
Suppose that
\[
G_{m-1}'(y_1)=G_{m-1}'(y_2)=\cdots=G_{m-1}'(y_m)=0,
\]
where $y_1<y_2<\cdots<y_m<0$.
By~\eqref{Gm}, the function $G_m(x)$ has $m+2$ roots
\[
-\infty,\ y_1,\ y_2,\ \ldots,\ y_m,\ 0.
\]
Because of~\eqref{G}, we see that $y_1,y_2,\ldots,y_m$
are $m$ distinct negative roots of
$M_m(x)$. This completes the proof. \qed

It should be mentioned that Theorem~\ref{thm_RootM}
can also be deduced from the criteria of
Liu and Wang~\cite{LW07}.
The following theorem gives an estimate of $M_n$.

\begin{thm}\label{thm_Mn}
We have
\begin{equation}\label{Mn}
M_n=\frac{1}{\sqrt{2r_1+1}}
\exp\biggl(2nr_1-n+\frac{n}{2r_1}+2r_1-1\biggr)
\biggl[1+O\Bigl(\frac{\log^{7/2}{n}}{\sqrt{n}}\Bigr)\biggr],
\end{equation}
where $r_1$ is the unique positive solution of the equation
$r_1\bigl(e^{2r_1}+1\bigr)=n$.
\end{thm}

\pf
Let $r=r_1$.
Applying Cauchy's formula and the generating function~\eqref{EgfM}, we have
\begin{equation}\label{Cauchy'}
\frac{M_n}{n!}=\frac{1}{2\pi i}\int_{|z|=r}\frac{F(z)}{z^{n+1}}\,dz
=\frac{1}{2\pi r^n\sqrt{e}}\,\int_{|\theta|\le\pi}e^Ad\theta,
\end{equation}
where
\begin{equation}\label{A'}
A=\frac{1}{2}e^{2r\cos{\theta}}+re^{i\theta}-n\theta\,i.
\end{equation}
We divide the integral in~\eqref{Cauchy'} into two parts as
\begin{equation}\label{2parts}
\int_{|\theta|\le\pi}e^Ad\theta=\int_{|\theta|\le\theta_0}e^Ad\theta+
\int_{\theta_0\le|\theta|\le\pi}e^Ad\theta,
\end{equation}
where
\[
\theta_0=\sqrt{{2\log{n}\over n}}.
\]
Let $\Re(A)$ denote the real part of $A$,
and $\Im(A)$ the imaginary part.
It follows from~\eqref{A'} that
\begin{align}
\Re(A)&=\frac{1}{2}e^{2r\cos{\theta}}\cos(2r\sin{\theta})
+r\cos{\theta},\label{ReA}\\
\Im(A)&=\frac{1}{2}e^{2r\cos{\theta}}
\sin(2r\sin{\theta})+r\sin{\theta}-n\theta.\notag
\end{align}
For the part $\int_{|\theta|\le\theta_0}\!e^Ad\theta$,
we have
\begin{align}
\Re(A)&={e^{2r}+2r\over 2}-{n(2r+1)\over2}\theta^2
+O\bigl(nr^3\theta_0^4\bigr),\label{re1}\\[5pt]
\Im(A)&=O\bigl(nr^2\theta_0^3\bigr).\notag
\end{align}
Substituting them into $e^A=e^{\Re(A)+i\,\Im(A)}$, we get
\begin{equation}\label{integral}
\int_{|\theta|\le\theta_0}e^Ad\theta
=\exp\Bigl({e^{2r}+2r\over 2}\Bigr)\int_{|\theta|\le\theta_0}e^{-m\theta^2}d\theta\,
\biggl(1+O\bigl(nr^2\theta_0^3\bigr)\biggr),
\end{equation}
where $m=(2r+1)n/2$. Note that
\[
\int_{x}^{\infty}e^{-t^2}dt=o\bigl(e^{-x^2}\bigr),\qquad
{\rm as}\quad x\to\infty.
\]
The integral in~\eqref{integral} can be estimated as follows
\begin{equation}\label{9}
\int_{|\theta|\le\theta_0}e^{-m\theta^2}d\theta
=\frac{1}{\sqrt{m}}\biggl(\sqrt{\pi}-2\int_{\sqrt{m/n}}^{\infty}\,e^{-t^2}dt\biggr)
=\sqrt{{\pi\over m}}\,\Bigl(1+o\bigl(e^{-r}\bigr)\Bigr).
\end{equation}
By~\eqref{integral} and~\eqref{9}, we find
\begin{equation}\label{47}
\int_{|\theta|\le\theta_0}e^Ad\theta
=\exp\Bigl({e^{2r}+2r\over 2}\Bigr)\sqrt{\frac{2\pi}{(2r+1)n}}\,\Bigl(1+O\bigl(nr^2\theta_0^3\bigr)\Bigr).
\end{equation}
Now we estimate the integration
$\int_{\theta_0\le|\theta|\le\pi}e^Ad\theta$.
By~\eqref{ReA}, we have
\[
\int_{\theta_0\le|\theta|\le\pi}e^Ad\theta
\le2\pi\max_{\theta_0\le\theta\le\pi}e^{\Re(A)}
\le2\pi\exp\Bigl(\frac{1}{2}e^{2r\cos{\theta_0}}+r\Bigr).
\]
Since
\[
2r\cos{\theta_0}=2r-r\theta_0^2+O\bigl(r\theta_0^4\bigr),
\]
we get
\[
\int_{\theta_0\le|\theta|\le\pi}e^Ad\theta
=O\biggl(\exp\Bigl({e^{2r}\over 2}-{n\theta_0^2\over 2}+r\Bigr)\biggr).
\]
It is easy to check that
\begin{equation}\label{6}
\lim_{n\to\infty}
\frac{\exp\bigl({e^{2r}\over 2}-{n\theta_0^2\over 2}+r\bigr)}
{\exp\bigl({e^{2r}+2r\over
2}\bigr)\sqrt{\frac{2\pi}{(2r+1)n}}\,nr^2\theta_0^3}
=0.
\end{equation}
Namely, the remainder of $\bigl|\int_{\theta_0\le|\theta|\le\pi}e^Ad\theta\bigr|$
is smaller than the remainder of
$\bigl|\int_{|\theta|\le\theta_0}e^Ad\theta\bigr|$.
By~\eqref{47}, we have
\[
\int_{|\theta|\le\pi}e^Ad\theta
=\exp\Bigl({e^{2r}+2r\over 2}\Bigr)\sqrt{\frac{2\pi}{(2r+1)n}}\,\Bigl(1+O\bigl(nr^2\theta_0^3\bigr)\Bigr).
\]
Hence by~\eqref{Cauchy'} and Stirling's formula
\[
n!=\frac{\sqrt{2\pi n}\,n^n}{e^n}\Bigl(1+O\bigl(n^{-1}\bigr)\Bigr),
\]
we have
\begin{equation}\label{1}
M_n=\frac{1}{\sqrt{2r+1}}\left({n\over r}\right)^n
\exp\Bigl(\frac{n}{2r}-n+r-1\Bigr)
\biggl[1+O\Bigl(\frac{\log^{7/2}{n}}{\sqrt{n}}\Bigr)\biggr].
\end{equation}
By Equation~\eqref{r1} and Lemma~\ref{lem_r}, we find
\[
\left({n\over r}\right)^n=e^{2nr+r}\biggl(1+O\Bigl({\log^2n\over n}\Bigr)\biggr).
\]
Together with~\eqref{1}, we arrive at~\eqref{Mn}.
This completes the proof. \qed

As will be seen in the next theorem,
the remainder $O\bigl(\frac{\log^{7/2}{n}}{\sqrt{n}}\bigr)$ plays an essential role in estimating the variance $V(\xi_n)$.

\begin{thm}\label{thm_EV_M}
We have
\begin{align}
\E(\xi_n)&={M_{n+1}\over 2M_n}-1\sim{n\over\log{n}},
\label{EM}\\[5pt]
\V(\xi_n)&={M_{n+2}\over 4M_n}-{M_{n+1}^2\over 4M_n^2}-{1\over 2}
\sim{n\over \log^2{n}}.
\label{VM}
\end{align}
\end{thm}

\pf
It can be easily checked that the expectation and
the variance of $\xi_n$ can be expressed by
\begin{align*}
\E(\xi_n)&={M_n'(1)\over M_n},\\
\V(\xi_n)&=\E(\xi_n)-\E(\xi_n)^2+\frac{M_n''(1)}{M_n}.
\end{align*}
Thus we can deduce the exact formulas in~\eqref{EM} and~\eqref{VM}.
In view of Theorem~\ref{thm_Mn},
Lemma~\ref{lem_r} and Lemma~\ref{lem1},
we find
\[
{M_{n+1}\over 2M_n}-1\sim{n\over\log{n}}.
\]

We now proceed to derive the approximation in~\eqref{VM}.
Suppose that
\[
t_i(e^{2t_i}+1)=n+i,
\]
for $i=0,1,2$. By Theorem~\ref{thm_Mn}, we have
\begin{equation}\label{75}
{M_{n+2}\over M_n}-{M_{n+1}^2\over M_n^2}
=\left(\sqrt{2t_0+1\over 2t_2+1}e^A-{2t_0+1\over 2t_1+1}e^B\right)
\biggl(1+O\Bigl({\log^{7/2}{n}\over\sqrt{n}}\Bigr)\biggr),
\end{equation}
where
\begin{align*}
A&=4t_2+\Bigl(2nt_2-2nt_0-2+{1\over t_2}\Bigr)-\Bigl({n\over 2t_0}-{n\over 2t_2}\Bigr)
+(2t_2-2t_0),\\[5pt]
B&=4t_1+\Bigl(4nt_1-4nt_0-2+{1\over t_1}\Bigr)-\Bigl({n\over t_0}-{n\over t_1}\Bigr)
+(4t_1-4t_0).
\end{align*}
By Lemma~\ref{lem_r},
both $\sqrt{2t_0+1\over 2t_2+1}$ and ${2t_0+1\over 2t_1+1}$
can be estimated by $1+O\bigl({1\over n\log{n}}\bigr)$.
Because of the estimates in Lemma~\ref{lem1},
\eqref{75} simplifies to
\begin{equation}\label{72}
{M_{n+2}\over M_n}-{M_{n+1}^2\over M_n^2}
=\bigl(e^A-e^B\bigr)\biggl(1+O\Bigl({\log^{7/2}{n}\over\sqrt{n}}\Bigr)\biggr).
\end{equation}
By Cauchy's mean value theorem, there exists a constant $C$ such that $B<C<A$ and
\begin{equation}\label{73}
e^A-e^B=(A-B)e^C.
\end{equation}
On one hand, Lemma~\ref{lem1} yields
\begin{align}
A-B&=\Bigl(4t_2-4t_1+{1\over t_2}-{1\over t_1}\Bigr)-(2n+2)(2t_1-t_0-t_2)
+{n\over 2}\Bigl({1\over t_0}+{1\over t_2}-{2\over t_1}\Bigr)\notag\\[5pt]
&={1\over n}\biggl(1+O\Bigl({1\over n^2\log{n}}\Bigr)\biggr).\label{74}
\end{align}
On the other hand, by Lemma~\ref{lem_r} we find that
\begin{equation}\label{e^C_M}
e^C={4n^2\over\log^2{n}}\biggl(1+O\Bigl({\log\log{n}\over\log{n}}\Bigr)\biggr).
\end{equation}
Substituting~\eqref{e^C_M} and~\eqref{74} into~\eqref{73}, we deduce that
\begin{equation}\label{71}
e^A-e^B={4n\over\log^2{n}}\biggl(1+O\Bigl({\log\log{n}\over\log{n}}\Bigr)\biggr).
\end{equation}
Substituting~\eqref{71} into~\eqref{72},
we obtain the approximation of $V(\xi_n)$.
This completes the proof.
\qed

By~\eqref{VM},
we see that $V(\xi_n)$ tends to infinity as $n\to\infty$.
Hence Theorem~\ref{thm_DistM} follows from Theorem~\ref{thm_RootM} and Proposition~\ref{prop_1}.

For $B_n$-partitions without zero-block,
we have an analogous limiting distribution.
Using the saddle point method as in the proof of Theorem~\ref{thm_Mn},
we obtain the following estimates of $N_n$.

\begin{thm}\label{thm_Nn}
We have
\begin{align}
N_n&=\frac{1}{\sqrt{(2r_0+1)}}\,\exp\biggl(2nr_0-n+\frac{n}{2r_0}-{1\over2}\biggr)
\biggl(1+O\Bigl(\frac{\log^{7/2}{n}}{\sqrt{n}}\Bigr)\biggr)\label{eqNn}\\[5pt]
&\sim{1\over\sqrt{\log{n}}}\,\exp\biggl(2nr_0-n+{n\over 2r_0}-{1\over2}\biggr).
\label{cor_Nn}
\end{align}
where $r_0$ is the unique positive solution of the equation
$r_0e^{2r_0}=n$.
\end{thm}

We remark that the approximation~\eqref{cor_Nn} can also be proved by
Hayman's theorem~\cite{Hay56}.

\begin{cor}
We have
\begin{equation}\label{N/M}
{N_n\over M_n}\sim\sqrt{\,{\log{n}\over 2n}}\,.
\end{equation}
\end{cor}

\pf Let $r_0e^{2r_0}=n$ and $r_1(e^{2r_1}+1)=n$.
By Theorem~\ref{thm_Mn} and Lemma~\ref{lem_r}, we obtain that
\[
M_n\sim{2n\over \log^{3/2}{n}}\,
\exp\biggl(2nr_1-n+\frac{n}{2r_1}-1\biggr).
\]
Using~\eqref{cor_Nn}, we get
\begin{equation}\label{2}
{N_n\over M_n}\sim\frac{\log{n}}{2n}
\exp\biggl(2n(r_0-r_1)-{n(r_0-r_1)\over2r_0r_1}+{1\over2}\biggr).
\end{equation}
By Cauchy's mean value theorem, we have
\begin{equation}\label{52}
n(r_0-r_1)={r_0\over 2}-{1\over 4}+O\Bigl({1\over\log{n}}\Bigr).
\end{equation}
Thus~\eqref{N/M} follows from~\eqref{2} and~\eqref{52}. This completes the proof. \qed

Recall that $N_{n,\,k}$ is the number of $B_n$-partitions without zero-block
having $k$ block pairs.
It can be verified that for any $n\geq 1$, the polynomial
\[
N_n(x)=\sum_{k}N_{n,\,k}x^k
\]
has $n$ distinct real roots.
Let $\xi_n'$ be the random variable of
the number of block pairs in $B_n$-partitions without zero-block.
Using the same argument as that for $\xi_n$,
we find
\begin{align*}
E(\xi_n')&={N_{n+1}\over 2N_n}-1
\sim{n\over\log{n}},\\[5pt]
V(\xi_n')&={N_{n+2}\over 4N_n}-{N_{n+1}^2\over 4N_n^2}-{1\over 2}
\sim{n\over \log^2{n}}.
\end{align*}
Hence $\V(\xi_n')$ tends to infinity as $n$ does.
By Proposition~\ref{prop_1}, we are led to the following assertion.

\begin{thm}\label{thm_DistN}
The limiting distribution of the random variable $\xi_n'$
is normal.
\end{thm}

\end{document}